\documentclass[10pt]{article}
\usepackage{amsmath}
\usepackage{amsfonts}
\usepackage{stmaryrd}
\usepackage{mathrsfs}
\usepackage{amssymb,amsfonts,amsmath, psfrag,eepic,colordvi,epsfig}
\usepackage{enumerate}
\usepackage{cite}
\usepackage{hyperref}
\topskip 
\numberwithin{equation}{section}
\newtheorem{theorem}{Theorem}[section]

\newtheorem{conjecture}[theorem]{Conjecture}

\newtheorem{lemma}[theorem]{Lemma}

\begin{document}
\parskip 7pt

\pagenumbering{arabic}
\def\sof{\hfill\rule{2mm}{2mm}}
\def\ls{\leq}
\def\gs{\geq}
\def\SS{\mathcal S}
\def\qq{{\bold q}}
\def\MM{\mathcal M}
\def\TT{\mathcal T}
\def\EE{\mathcal E}
\def\lsp{\mbox{lsp}}
\def\rsp{\mbox{rsp}}
\def\pf{\noindent {\it Proof.} }
\def\mp{\mbox{pyramid}}
\def\mb{\mbox{block}}
\def\mc{\mbox{cross}}
\def\qed{\hfill \rule{4pt}{7pt}}
\def\pf{\noindent {\it Proof.} }
\textheight=22cm

\begin{center}
{\Large\bf   Proof of a conjecture of Keith
 on congruences
  of the reciprocal  of a
  false theta function  }
\end{center}

\begin{center}

Jing Jin$^{1}$,   Sijia Wang$^{2}$
  and Olivia X.M. Yao$^{3}$

  $^{1}$College of Agricultural Information,\\
   Jiangsu Agri-animal
  Husbandry Vocational College,
  \\
  Taizhou, 225300,  Jiangsu,  P. R. China

$^{2,3}$School of Mathematical Sciences, \\
  Suzhou University of Science and
Technology, \\
 Suzhou,  215009, Jiangsu Province,
 P. R. China

email: $^{1}$jinjing19841@126.com,
 $^{2}$Jia\_hao1109@163.com,
   $^{3}$yaoxiangmei@163.com

 \end{center}

\noindent {\bf Abstract.}
 Recently, Keith investigated
  reciprocals of false theta functions
   and proved some interesting results
   such as congruences,  asymptotic bounds,
    and combinatorial identities.
     At the end of his paper, Keith posed
      a conjecture on congruences
       modulo 4 and 8 for the coefficients
        of the  reciprocal  of a
         false theta function.
          In this paper, we
           not only confirm Keith's
            conjecture, but also prove
             a generalized result.



\section{Introduction}

\allowdisplaybreaks

Ramanujan's
  general  theta
     function is  defined by
\begin{align}\label{1-1}
f(a,b)=\sum_{n=-\infty}^\infty
 a^{n(n+1)/2}b^{n(n-1)/2}, \qquad |ab|<1.
\end{align}
In Ramanujan's notation,
 the Jacobi triple product
 identity takes the shape
\begin{align}\label{1-2}
f(a,b)=(-a,-b,ab;ab)_\infty.
 \end{align}
Here and throughout this paper,
 we use the following
  notation:
\begin{align*}
 (a;q)_\infty &\;:=\prod_{k=0}^\infty
 (1-aq^k), \qquad
 (a;q)_n  :=\frac{(a;q)_\infty}
 {(aq^n;q)_\infty},\\
 (a_1,a_2,\ldots,a_k;q)_\infty
 &\;:=(a_1;q)_\infty (a_{2};q)_\infty
 \cdots (a_k;q)_\infty,
\end{align*}
where $q$ is a complex number
 with $|q|<1$. In addition, for all positive
  integers $m$, define
  \begin{align*}
 f_m:= (q^m;q^m)_\infty.
  \end{align*}

Three special  cases
  of \eqref{1-1} are defined  by, in Ramanujan's
   notation,
\begin{align}
f( -q,-q^2)&\;=\sum_{n=-\infty}^\infty
 (-1)^n q^{n(3n-1)/2}=f_1,\label{1-3}\\
f(-q,-q)&\;=\sum_{n=-\infty}^\infty
 (-1)^n q^{n^2}=\frac{f_1^2}{f_2},\label{1-4}\\
 f(-q,-q^3)&\;=\sum_{n=0}^\infty
 (-q)^{n(n+1)/2}
  =\frac{f_1f_4}{f_2}. \label{1-5}
\end{align}
The reciprocals of the right
 sides
 of \eqref{1-3}--\eqref{1-5}
    are the generating functions of
 the ordinary partition function $p(n)$ \cite{Andrews-1976},
  the  overpartition
  function $\overline{p}
  (n)$
    \cite{Corteel},
 and the partition
 with odd parts distinct
 function \cite{Hirschhorn-1}, respectively,
  which are three of the most important
types of partition functions. A number of
  congruences for the three types of
partition functions have been proved;
 see for example
 \cite{HS-2005b,Hirschhorn-1,Ramanujan-1}.

In 1917, Rogers \cite{Rogers} introduced false theta functions,
    which are series
that would be classical theta functions
 except for changes in signs
of an infinite number of terms. False theta functions closely
resemble ordinary theta functions; however, they do not have the
modular transformation properties that theta functions have.
  In his  notebooks \cite{Ramanujan-2}
   and lost notebook \cite{Ramanujan-3}, Ramanujan
  recorded many false theta function
   identities on false theta functions
    $\Psi(a,b)$ defined by
\begin{align}\label{1-6}
\Psi(a,b):=\sum_{n=0}^\infty a^{n(n+1)/2} b^{n(n-1)/2}
-\sum_{n=-\infty}^{-1}
 a^{n(n+1)/2}b^{n(n-1)/2}.
\end{align}

  In recent times, there has been an increase in
interest in false theta functions. They appear in a variety of
contexts including modular and quantum modular forms,
 meromorphic Jacobi forms,
 quantum knot invariants,
 and many more.
 In \cite{Andrews-1981},   Andrews   proved five identities
related to   false theta functions
 which were recorded in Ramanujan's
   lost notebook \cite{Ramanujan-3} (c.f.
  \cite[Sect. 9.3, pp. 227--232]{Andrews-2005}).  Other
analytic proofs have been given by Andrews
 and
Warnaar
 \cite{Andrews-2007}
 and by Wang\cite{Wang}.
 Burson \cite{Burson}
gave a bijection proof
 of the  following identity due to Ramanujan:
  \[
\Psi(q^3,q)=\sum_{n=0}^\infty \frac{ (q;q^2)_n
 }{(-q;q^2)_{n+1}}q^n.
  \]

Recently, Keith \cite{Keith} investigated
   reciprocals
  of false theta functions
   and proved some interesting results
    on congruences,
     asymptotic bounds,
     and combinatorial identities.
    For example,  he proved
       that for $n\geq 0$,
\begin{align*}
c_5(8n+5) &\; \equiv 0 \pmod 2,  \\
c_5(32n+31) &\; \equiv 0 \pmod 4,\\
c_9(16n+12)&\;\equiv 0 \pmod 2,
\end{align*}
where $c_t(n)$ is defined by
\begin{align}\label{1-7}
\sum_{n=0}^\infty c_t(n)q^n:=\frac{1}{\Psi(-q^t,q)}.
\end{align}

At the end of his paper, Keith \cite{Keith}
 posed the following conjecture
  on congruences modulo 4 and 8 for $c_5(n)$:

 \begin{conjecture} \cite{Keith}
 \label{conj-1}
For $n\geq 0$,
\begin{align}
c_5(32n+31) &\;\equiv 0 \pmod 8,\label{1-8}\\[6pt]
c_5(128n+123) &\;\equiv 0 \pmod 8,\label{1-9}\\[6pt]
c_5(512n+491) &\;\equiv 0 \pmod 8,\label{1-10}\\[6pt]
c_5(64n+19) &\;\equiv 0 \pmod 4,\label{1-11}\\[6pt]
c_5(256n+75) &\;\equiv 0 \pmod 4,\label{1-12}\\[6pt]
c_5(196n+i) &\;\equiv 0 \pmod 4,\label{1-13}\\[6pt]
c_5(196n+j) &\;\equiv 0 \pmod 4, \label{1-14}
\end{align}
where $i\in \{110,138,194\}$
 and $j\in\{19,47,75,
  103,159,187\}$.
\end{conjecture}

The aim of this paper is
 to
  prove congruences \eqref{1-8}--\eqref{1-13}
    and generalize \eqref{1-14}.
 The main results
  of this paper can be stated   as follows.

\begin{theorem}\label{Th-1}
Congruences \eqref{1-8}--\eqref{1-13} are true.
\end{theorem}

\begin{theorem}\label{Th-2}
Let $p\equiv 7\pmod 8 $ be a prime. For    all   integers
  $n$ and $k\geq 0$ with $p\nmid n$,
\begin{align}\label{1-15}
c_5\left(
 4p^{2k+1}n+\frac{8p^{2k+2}+1}{3}
\right)\equiv 0 \pmod 4.
\end{align}
\end{theorem}

Setting $p=7$ in \eqref{1-15}
  yields \eqref{1-14}. Therefore, \eqref{1-15} is a
generalization
 of \eqref{1-14}.

 This paper is organized as follows.
               Sections  2 and 3 are  devoted to
         the proofs of Theorems
          \ref{Th-1} and \ref{Th-2},
          respectively.
          In Section 4, we make some
          concluding remarks
          concerning future
directions.

\section{Proof of Theorem \ref{Th-1}}

The aim of this section  is to pose a
 proof of Theorem \ref{Th-1}.
Throughout this paper,
 we always set
\begin{align}\label{2-0}
A:=\sum_{n=0}^\infty (-1)^{n(n+1)/2}
 q^{3n^2+2n},\qquad
 B:=\sum_{n=-\infty}^{-1} (-1)^{n(n+1)/2}
  q^{3n^2+2n}.
\end{align}

We first prove some lemmas.

\begin{lemma}\label{L-1}
Define
\begin{align}\label{2-1}
\sum_{n=0}^\infty
 b_1(n)q^n:=\frac{1}{A+B}.
\end{align}
Then for $n\geq 0$,
\begin{align}
b_1(32n+31)\equiv&\; 0 \pmod 8,\label{2-2}\\
b_1(128n+123)\equiv&\; 0 \pmod 8,\label{2-3}\\
b_1(512n+491)\equiv&\; 0 \pmod 8,\label{2-4}\\
b_1(64n+19)\equiv&\; 0 \pmod 4,\label{2-5}\\
b_1(256n+75)\equiv&\; 0 \pmod 4.\label{2-6}
\end{align}
\end{lemma}

\noindent{\it Proof.}  In view of \eqref{1-1},
  \eqref{1-2}
 and \eqref{2-0},
 \begin{align*}
A+B= &\; f(-q^5,q)
 = (-q,q^5,-q^6;-q^6)_\infty
 \nonumber\\
 =&\;(-q,q^5,q^7, -q^{11};q^{12})_\infty
 \frac{f_{12}^3}{f_6f_{24}}.
 \end{align*}
Therefore,
\begin{align}\label{2-8}
\sum_{n=0}^\infty
 b_1(n)q^n=&\;\frac{f_6f_{24}}{
 (-q,q^5,q^7,-q^{11};q^{12})_\infty
  f_{12}^3}
  \nonumber\\
  =&\;\frac{f_6f_{24} (q,-q^5,-q^7,q^{11};q^{12})_\infty
   }{
 ( q^2,q^{10},q^{14},q^{22};q^{24})_\infty
  f_{12}^3}
  \nonumber\\
  =&\;\frac{f_4f_6^2f_{24}}{
   f_2f_{12}^6}f(-q,-q^{11})f(q^5,q^7).
\end{align}
It follows from \cite[Entry 29, pp. 45]{Berndt-1} that
 if $ab=cd$, then
 \begin{align}\label{2-9}
f(a,b)f(c,d)=f(ac,bd)f(ad,bc) +af(b/c,ac^2d)
 f(b/d,acd^2).
 \end{align}
Setting $a=-q,b=-q^{11},c=q^5,d=q^7$
 into \eqref{2-9} yields
 \begin{align}\label{2-10}
f(-q,-q^{11})f(q^5,q^7)= \frac{ f_6f_8f_{24}
 }{f_{12}}-q\frac{f_4f_6f_{24}^3
  }{f_8f_{12}^2}.
 \end{align}
Thanks to  \eqref{2-8} and \eqref{2-10},
\begin{align}
\sum_{n=0}^\infty
 b_1(n)q^n=\frac{1}{A+B}
  =&\;\frac{f_4f_6^3f_{24}^2}{
   f_2f_{12}^7}\left(
    f_8-q\frac{f_4f_{24}^2
    }{f_8f_{12}}\right),\label{2-11-0}
\end{align}
from which, we arrive at
\begin{align}
\sum_{n=0}^\infty b_1(2n+1)q^n=&\; -\frac{f_3^3}{f_1} \cdot
\frac{f_2^2
f_{12}^4}{f_4f_6^8}\nonumber\\
\equiv &\; -\frac{f_3^3}{f_1} \cdot
 \frac{f_2^2  }{f_4 } \pmod 8.\label{2-11}
\end{align}
Here we have used the fact that
 for all positive
  integers $m$ and $k$
\begin{align}\label{2-12}
f_k^{2^m}\equiv f_{2k}^{2^{m-1}} \pmod {2^m}.
 \end{align}
Xia and Yao \cite[(2.15)]{Yao-1} proved
 that
\begin{align}\label{2-13}
\frac{f_3^3}{f_1}=\frac{ f_4^3
 f_6^2}{f_2^2f_{12}}+q\frac{f_{12}^3
  }{f_4}.
 \end{align}
Substituting \eqref{2-13}
  into \eqref{2-11} and
  extracting those
terms in which the power of  $q
 $ is congruent to 1 modulo 2, then
 dividing  them by  $q$  and replacing
$q^{2}$ by $q$ yields
 \begin{align}\label{2-15}
\sum_{n=0}^\infty b_1(4n+3)q^n\equiv
 -f_1^2\cdot \frac{f_6^3}{f_2^2} \pmod 8.
 \end{align}
It follows from
\cite[Corollary (i) and
 (ii), pp. 49]{Berndt-1}
that
 \begin{align}\label{2-16}
f_1^2=\frac{f_2f_8^5}{f_4^2f_{16}^2
 }-2q\frac{f_2f_{16}^2}{f_8}
 \end{align}
 and
 \begin{align}\label{2-17}
\frac{1}{f_1^2}= \frac{f_8^5}{ f_2^5f_{16}^2
 }+2q\frac{f_4^2f_{16}^2
  }{f_2^5f_8}.
 \end{align}
If we substitute \eqref{2-16}
  into \eqref{2-15}, we arrive at
\[
\sum_{n=0}^\infty b_1(4n+3)q^n\equiv
 -\frac{f_6^3f_8^5}{
 f_2f_4^2f_{16}^2}+2q\frac{f_6^3f_{16}^2
 }{f_2f_8} \pmod 8,
\]
 which yields
\begin{align}\label{2-18}
\sum_{n=0}^\infty b_1(8n+3)q^n\equiv
 -\frac{f_3^3}{f_1}
 \cdot
 \frac{f_4^5}{
 f_2^2f_{8}^2} \pmod 8
\end{align}
and
\begin{align}\label{2-19}
\sum_{n=0}^\infty b_1(8n+7)q^n\equiv
 2\frac{f_3^3}{f_1}\cdot
   \frac{f_8^2}{f_4} \pmod 8.
 \end{align}
Substituting \eqref{2-13}
  into \eqref{2-19} and
  extracting those
terms in which the power of  $q
 $ is congruent to 1 modulo 2, then
 dividing  them by  $q$  and replacing
$q^{2}$ by $q$, we obtain
  \begin{align*}
\sum_{n=0}^\infty b_1(16n+15)q^n\equiv
 2
   \frac{f_4^2f_6^3}{f_2^2} \pmod 8,
 \end{align*}
which yields \eqref{2-2}.

Substituting \eqref{2-13}
  into \eqref{2-18} yields
\[
\sum_{n=0}^\infty b_1(8n+3)q^n\equiv
 -\frac{f_4^8f_6^2}{f_2^4f_8^2f_{12}}
 -q\frac{f_4^4f_{12}^3
  }{f_2^2 f_8^2} \pmod 8,
\]
from which, we arrive at
\begin{align}
\sum_{n=0}^\infty b_1(16n+3)q^n\equiv -
\frac{f_2^8f_3^2}{f_1^4f_4^2f_{6}} \pmod 8
\label{2-20}
\end{align}
and
\begin{align}
\sum_{n=0}^\infty b_1(16n+11)q^n\equiv
 -\frac{1}{f_1^2}\cdot
  \frac{f_2^4f_6^3}{
  f_{4}^2} \pmod 8.\label{2-21}
\end{align}
If we substitute \eqref{2-17} into \eqref{2-21}, we get
\begin{align*}
\sum_{n=0}^\infty b_1(16n+11)q^n\equiv
 -\frac{f_6^3f_8^5}{
 f_2f_4^2f_{16}^2}-2q\frac{f_6^3f_{16}^2
  }{f_2f_8} \pmod 8,
\end{align*}
which yields
\begin{align}\label{2-22}
\sum_{n=0}^\infty b_1(32n+11)q^n\equiv
 -\frac{f_3^3f_4^5}{
 f_1f_2^2f_{8}^2} \pmod 8
\end{align}
and
\begin{align}
\sum_{n=0}^\infty b_1(32n+27)q^n\equiv
 -2\frac{f_3^3f_8^2}{
 f_1f_{4}} \pmod 8. \label{2-23}
\end{align}
Combining \eqref{2-19} and \eqref{2-23} yields
\[
b_1(32n+27)\equiv -b_1(8n+7) \pmod 8.
\]
Replacing $n$ by $4n+3$
 in the above congruence
   and using \eqref{2-2}, we get \eqref{2-3}.

In view of \eqref{2-18} and \eqref{2-22},
\begin{align}\label{a-1}
b_1(32n+11)\equiv b_1(8n+3) \pmod 8.
\end{align}
Replacing $n$ by $16n+15$ in
 the above congruence
   and using \eqref{2-3},
 we get \eqref{2-4}.

 In view of \eqref{2-12} and \eqref{2-20},
 \begin{align}
\sum_{n=0}^\infty b_1(16n+3)q^n\equiv - f_3^2
\cdot
 \frac{f_2^2 }{ f_{6}} \pmod 4. \label{2-24}
\end{align}
Substituting \eqref{2-16}
  into \eqref{2-24} and
  extracting those
terms in which the power of  $q
 $ is congruent to 1 modulo 2, then
 dividing  them by  $q$  and replacing
$q^{2}$ by $q$, we obtain
 \begin{align*}
\sum_{n=0}^\infty b_1(32n+19)q^n\equiv&\;  2q\frac{f_1^2f_{24}^2
 }{f_{12}}\nonumber\\
 \equiv &\; 2qf_2\frac{f_{24}^2}{f_{12}}\pmod
4. \qquad ({\rm by}\ \eqref{2-12})
\end{align*}
Congruence  \eqref{2-5} follows from the above congruence.

Replacing $n$ by $8n+2$ in \eqref{a-1}
  and utilizing \eqref{2-5},
 we get \eqref{2-6}. This completes
  the proof
   of Lemma \ref{L-1}. \qed

   \begin{lemma}\label{L-2}
Define
\begin{align}\label{2-26}
\sum_{n=0}^\infty
 b_2(n)q^n:=\frac{B}{(A+B)^2}.
\end{align}
Then for $n\geq 0$,
\begin{align}
b_2(32n+31)\equiv&\; 0 \pmod 4,
\label{2-27}\\
b_2(128n+123)\equiv&\; 0 \pmod 4,\label{2-28}\\
b_2(512n+491)\equiv&\; 0 \pmod 4,\label{2-29}\\
b_2(32n+19)\equiv&\; 0 \pmod 2,\label{2-30}\\
b_2(128n+75)\equiv&\; 0 \pmod 2. \label{2-31}
\end{align}

\end{lemma}

\noindent{\it Proof.}
 In light of \eqref{2-0},
  \eqref{2-11-0} and \eqref{2-26},
\begin{align}\label{2-32}
\sum_{n=0}^\infty
 b_2(n)q^n&\; =\frac{f_4^2f_6^6f_{24}^4}{
   f_2^2f_{12}^{14}}\left(
    f_8-q\frac{f_4f_{24}^2
    }{f_8f_{12}}\right)^2\sum_{n=1}^\infty
     (-1)^{n(n-1)/2}q^{3n^2-2n}
     \nonumber\\
&\; \equiv \frac{f_2^2  }{
   f_6^6}\left(
    f_8-q\frac{f_4f_{24}^2
    }{f_8f_{12}}\right)^2 \nonumber\\
&\quad \times    \left(\sum_{n=1}^\infty
     (-1)^n q^{12n^2-4n}+q\sum_{n=0}^\infty
      (-1)^nq^{12n^2+8n}\right) \pmod 4.
      \qquad ({\rm by}\ \eqref{2-12})
\end{align}
Extracting those terms in which the power of  $q
 $ is congruent to 1 modulo 2 in \eqref{2-32}, then
 dividing  them by  $q$
  and replacing
$q^{2}$ by $q$, we obtain
\begin{align}
\sum_{n=0}^\infty
 b_2(2n+1)q^n
 \equiv &\;   q\frac{f_1^2}{f_3^2} \cdot
 \frac{  f_{12}^2}{
 f_2^2}\sum_{n=0}^\infty
  (-1)^n q^{6n^2+4n}+2f_4\sum_{n=1}^\infty
   (-1)^n q^{6n^2-2n}\nonumber\\
   &\; +\frac{f_1^2}{
   f_3^2} \cdot
   \frac{f_4^2}{f_6^2}
   \sum_{n=0}^\infty (-1)^n q^{6n^2+4n}
   \pmod 4.
   \qquad ({\rm by}\ \eqref{2-12}) \label{2-33}
\end{align}
Moreover,
  it follows from
   \cite[(3.29)]{Xia-Yao-2} that
  \begin{align}
  \label{a-2}
\frac{f_3^2}{f_1^2}&\;=
 \frac{f_4^4f_6f_{12}^2}{f_2^5
f_8f_{24}}+2q\frac{f_4f_6^2f_8f_{24}} {f_2^4f_{12}} .
\end{align}
Replacing $q$ by $-q$ in \eqref{a-2} yields
\begin{align}
\frac{f_1^2}{f_3^2}= \frac{ f_2f_4^2f_{12}^4
 }{f_6^5f_8f_{24}}
 -2q\frac{f_2^2f_8f_{12}f_{24}}{
 f_4f_6^4}.\label{2-34}
\end{align}
Substituting \eqref{2-34}
  into \eqref{2-33} and
  extracting those
terms in which the power of  $q
 $ is congruent to 1 modulo 2, then
 dividing  them by  $q$  and replacing
$q^{2}$ by $q$, we arrive at
\begin{align}
\sum_{n=0}^\infty
 b_2(4n+3)q^n
&\;\equiv   \left(2f_8+\frac{f_1^3}{f_3} \cdot \frac{f_{12
}}{f_4}\right) \left(\sum_{n=0}^\infty
  q^{12n^2+4n} -\sum_{n=0}^\infty
   q^{12n^2+16n+5}\right)
   \pmod 4. \qquad ({\rm by}\ \eqref{2-12}) \label{2-35}
\end{align}
In view of \cite[(2.14)]{Yao-1},
\begin{align}\label{2-36}
\frac{f_1^3}{f_3}= \frac{f_4^3}{f_{12}} -3q \frac{ f_2^2f_{12}^3
 }{f_4f_6^2}.
\end{align}
Substituting \eqref{2-36}
  into \eqref{2-35} yields
\begin{align*}
\sum_{n=0}^\infty
 b_2(4n+3)q^n
 \equiv &\;
    \left(2f_8+ f_4^2\right)
    \sum_{n=0}^\infty q^{12n^2+4n}
   -3q \frac{f_2^2f_{12}^4
     }{f_4^2f_6^2}
     \sum_{n=0}^\infty q^{12n^2+4n}
     \nonumber\\
     &\;-q^5(2f_8+f_4^2)  \sum_{n=0}^\infty
      q^{12n^2+16n}+3q^6
     \frac{f_2^2f_{12}^4
     }{f_4^2f_6^2}\sum_{n=0}^\infty
      q^{12n^2+16n}
  \pmod 4,
\end{align*}
from which, we obtain
\begin{align}
\sum_{n=0}^\infty
 b_2(8n+3)q^n
 \equiv &\;
    \left(2f_4+f_2^2\right)
    \sum_{n=0}^\infty q^{6n^2+2n}
+3q^3
     \frac{f_1^2}{f_3^2}
     \cdot\frac{f_{6}^4
     }{f_2^2 }\sum_{n=0}^\infty
      q^{6n^2+8n}
  \pmod 4 \label{2-37}
\end{align}
and
\begin{align}
\sum_{n=0}^\infty
 b_2(8n+7)q^n
&\;\equiv  -(q^2f_2^2+2q^2f_4) \sum_{n=0}^\infty
 q^{6n^2+8n}-3\frac{f_1^2}{f_3^2}
 \cdot \frac{f_6^4}{f_2^2
   }\sum_{n=0}^\infty
   q^{6n^2+2n}
 \pmod 4. \label{2-38}
\end{align}
If we substitute \eqref{2-34}
  into \eqref{2-38} and
  extract  those
terms in which the power of  $q
 $ is congruent to 1 modulo 2, then
 divide  them by  $q$  and replace
$q^{2}$ by $q$, we deduce that
\begin{align*}
\sum_{n=0}^\infty
 b_2(16n+15)q^n
&\;\equiv   2\frac{f_4f_6f_{12}}{f_2} \sum_{n=0}^\infty
 q^{3n^2+n}
   \pmod 4.
\end{align*}
Congruence \eqref{2-27}
  follows from the above congruence  since
 $3n^2+n$ is an even number.

Substituting \eqref{2-34}
  into \eqref{2-37} and
  extracting those
terms in which the power of  $q
 $ is congruent to 1 modulo 2, then
 dividing  them by  $q$  and replacing
$q^{2}$ by $q$, we get
\begin{align}
\sum_{n=0}^\infty
 b_2(16n+11)q^n
 \equiv &\;-
 \frac{q f_3^3 }{f_1 }
     \cdot\frac{f_2^2f_6^2
     }{ f_4f_{12} }\sum_{n=0}^\infty
      q^{3n^2+4n}
  \pmod 4. \qquad ({\rm by}\ \eqref{2-12}) \label{2-39}
\end{align}
Note that
\begin{align}\label{2-40}
\sum_{n=0}^\infty q^{3n^2+4n}=\sum_{n=0}^\infty
 q^{12n^2+8n}+q^7\sum_{n=0}^\infty
  q^{12n^2+20n}
\end{align}
and
\begin{align}\label{a-3}
\sum_{n=0}^\infty q^{3n^2+2n}=\sum_{n=0}^\infty
 q^{12n^2+4n}+q^5\sum_{n=0}^\infty
  q^{12n^2+16n}.
\end{align}
 Substituting \eqref{2-13} and \eqref{2-40}
  into \eqref{2-39} yields
\begin{align*}
\sum_{n=0}^\infty
 b_2(16n+11)q^n\equiv &\;-
  q\frac{f_4^2f_6^4}{f_{12}^2}
  \sum_{n=0}^\infty q^{12n^2+8n}
  -q^2\frac{f_2^2f_6^2f_{12}^2
   }{f_4^2}\sum_{n=0}^\infty q^{12n^2+8n}
   \nonumber\\
   &\;
   -q^8\frac{f_4^2f_6^4}{f_{12}^2}
    \sum_{n=0}^\infty q^{12n^2+20n}
    -q^9\frac{f_2^2f_6^2f_{12}^2
   }{f_4^2}\sum_{n=0}^\infty q^{12n^2+20n}
   \pmod 4,
\end{align*}
from which, we get
\begin{align}\label{2-41}
\sum_{n=0}^\infty
 b_2(32n+11)q^n
 \equiv &\;-
q\frac{f_1^2f_3^2 f_6^2}{f_2^2 }
 \sum_{n=0}^\infty q^{6n^2+4n}
 -q^4\frac{f_2^2f_3^4}{f_6^2}\sum_{n=0}^\infty
  q^{6n^2+10n}\nonumber\\
  \equiv &\;-
q\frac{ f_3^2 }{f_1^2 }\cdot f_6^2
 \sum_{n=0}^\infty q^{6n^2+4n}
 -q^4 f_2^2 \sum_{n=0}^\infty
  q^{6n^2+10n} \pmod 4 \qquad ({\rm by}\ \eqref{2-12})
\end{align}
and
\begin{align}\label{2-42}
\sum_{n=0}^\infty
 b_2(32n+27)q^n
 \equiv &\;-
\frac{f_2^2f_3^4}{f_6^2}\sum_{n=0}^\infty
 q^{6n^2+4n}-q^4\frac{f_1^2f_3^2f_6^2}{f_2^2}
 \sum_{n=0}^\infty q^{6n^2+10n}
 \nonumber\\
  \equiv &\;-
 f_2^2 \sum_{n=0}^\infty
 q^{6n^2+4n}-q^4\cdot
 \frac{f_3^2}{f_1^2}
 \cdot
  f_6^2
 \sum_{n=0}^\infty q^{6n^2+10n}
\pmod 4. \qquad ({\rm by}\ \eqref{2-12})
\end{align}
 Substituting \eqref{a-2}
  into \eqref{2-42} and
  extracting those
terms in which the power of  $q
 $ is congruent to 1 modulo 2, then
 dividing  them by  $q$  and replacing
$q^{2}$ by $q$, we arrive at
\begin{align*}
\sum_{n=0}^\infty
 b_2(64n+59)q^n
 \equiv &\;
2q^2 \frac{f_2f_3^4f_4f_{12}}{ f_1^4f_6}\sum_{n=0}^\infty
q^{3n^2+5n}
 \nonumber\\
  \equiv &\;
2q^2   \frac{ f_2 f_{12}^2}{
f_6}\sum_{n=0}^\infty q^{3n^2+5n}
\pmod 4.  \qquad ({\rm by}\ \eqref{2-12})
\end{align*}
Congruence \eqref{2-28} follows from the above congruence since
$3n^2+5n$
 is an even number.

If we substitute \eqref{a-2}
  into \eqref{2-41} and
  extract  those
terms in which the power of  $q
 $ is congruent to 1 modulo 2, then
 divide  them by  $q$  and replace
$q^{2}$ by $q$, we arrive at
\begin{align}
\sum_{n=0}^\infty
 b_2(64n+43)q^n
 \equiv &\;-
 \frac{f_3^3}{f_1}
 \cdot
  \frac{f_2^2f_6^2}{f_4f_{12}}
  \sum_{n=0}^\infty q^{3n^2+2n} \pmod 4.
  \qquad ({\rm by}\ \eqref{2-12}) \label{2-46}
\end{align}
Combining \eqref{2-39} and \eqref{2-46} yields
\begin{align}\label{2-47}
\sum_{n=0}^\infty
 b_2(16n+11)q^n+\sum_{n=0}^\infty
 b_2(64n+43)q^n
 \equiv &\;-
 \frac{  f_3^3 }{f_1 }
     \cdot\frac{f_2^2f_6^2
     }{ f_4f_{12} }\sum_{n=-\infty}^\infty
      q^{3n^2+2n}
  \pmod 4.
\end{align}
By \eqref{1-1}
 and \eqref{1-2},
\[
\sum_{n=-\infty}^\infty
      q^{3n^2+2n}=\frac{f_2^2
      f_3f_{12}}{f_1f_4f_6},
\]
from which with \eqref{2-12}
 and \eqref{2-47},
\begin{align}\label{2-48}
\sum_{n=0}^\infty
 b_2(16n+11)q^n+\sum_{n=0}^\infty
 b_2(64n+43)q^n
 \equiv &\;-
 \frac{f_2^4f_3^4f_6}{f_1^2f_4^2}
 \equiv-\frac{f_6^3}{f_1^2} \pmod 4.
\end{align}
 Substituting \eqref{2-17}
  into \eqref{2-48} and
  extracting those
terms in which the power of  $q
 $ is congruent to 1 modulo 2, then
 dividing  them by  $q$  and replacing
$q^{2}$ by $q$, we obtain
\begin{align}
\sum_{n=0}^\infty
 b_2(32n+27)q^n+\sum_{n=0}^\infty
 b_2(128n+107)q^n
 \equiv  2\frac{
 f_2^2f_3^3f_8^2}{f_1^5 f_4}
  \equiv 2\frac{f_3^3}{f_1}
  \cdot
   \frac{f_8^2}{f_4} \pmod 4.
   \qquad ({\rm by}\ \eqref{2-12}) \label{2-49}
\end{align}
Substituting \eqref{2-13}
  into \eqref{2-49} and
  extracting those
terms in which the power of  $q
 $ is congruent to 1 modulo 2, then
 dividing  them by  $q$  and replacing
$q^{2}$ by $q$ yields
\begin{align*}
\sum_{n=0}^\infty
 b_2(64n+59)q^n+\sum_{n=0}^\infty
 b_2(256n+235)q^n
 \equiv    2\frac{f_4^2 f_6^3}{f_2^2}
    \pmod 4,
\end{align*}
from which, we get
\[
 b_2(128n+123) +
 b_2(512n+491)\equiv 0 \pmod 4.
\]
Congruence \eqref{2-29} follows from \eqref{2-28}
  and the above congruence.

In view of \eqref{2-12} and \eqref{2-37},
\begin{align*}
\sum_{n=0}^\infty
 b_2(8n+3)q^n
 \equiv &\;
  f_4
    \sum_{n=0}^\infty q^{6n^2+2n}
+q^3 \frac{f_{6}^3
     }{f_2  }\sum_{n=0}^\infty
      q^{6n^2+8n}
  \pmod 2,
\end{align*}
which yields
\begin{align*}
\sum_{n=0}^\infty
 b_2(16n+3)q^n
 \equiv &\;
  f_2
    \sum_{n=0}^\infty q^{3n^2+ n}
  \pmod 2.
\end{align*}
Congruence \eqref{2-30} follows from the above congruence.

In light of \eqref{2-12} and \eqref{2-41},
\begin{align*}
\sum_{n=0}^\infty
 b_2(32n+11)q^n
 \equiv &\;
q\frac{   f_6^3}{f_2 }
 \sum_{n=0}^\infty q^{6n^2+4n}
 +q^4 f_4 \sum_{n=0}^\infty
  q^{6n^2+10n} \pmod 2,
\end{align*}
which yields
\begin{align*}
\sum_{n=0}^\infty
 b_2(64n+11)q^n
 \equiv &\;
 q^2  f_2 \sum_{n=0}^\infty
  q^{3n^2+5n} \pmod 2.
\end{align*}
Congruence \eqref{2-31}
  follows from  the above congruence
 since $3n^2+5n$ is an even number. This completes
  the proof of Lemma \ref{L-2}. \qed

\begin{lemma} \label{L-3}
Define
\begin{align}\label{2-51}
\sum_{n=0}^\infty
 b_3(n)q^n:=\frac{A^2B+AB^2}{(A+B)^4}.
\end{align}
Then for $n\geq 0$,
\begin{align}
b_3(16n+15)\equiv &\; 0 \pmod 2,\label{2-52}\\
b_3(64n+59)\equiv &\; 0 \pmod 2,\label{2-53}\\
b_3(256n+235)\equiv &\; 0 \pmod 2. \label{2-54}
\end{align}
\end{lemma}

\noindent{\it Proof}. By \eqref{1-1},
 \eqref{1-2} and \eqref{2-0},
 \begin{align}\label{2-55}
A+B=\sum_{n=-\infty}^\infty
 (-1)^{n(n+1)/2}q^{3n^2+2n}
  =(-q,q^5,-q^6;-q^6)_\infty\equiv (q,q^5,
q^6;q^6)_\infty \pmod 2.
 \end{align}
In addition, it is easy to check that
 \begin{align}\label{2-56}
(q,q^5,q^6;q^6)_\infty=\frac{f_1f_6^2}{f_2f_3}.
\end{align}
Combining \eqref{2-55} and \eqref{2-56} yields
 \begin{align}\label{2-57}
A+B\equiv \frac{f_1f_6^2}{f_2f_3} \pmod 2.
\end{align}
Thanks to \eqref{2-0} and \eqref{2-57},
\begin{align}
A^2B+AB^2\equiv&\; A^2\left(
 \frac{f_1f_6^2}{f_2f_3} -A  \right) +A
 \left( \frac{f_1f_6^2}{f_2f_3}
  -A  \right)^2
  \nonumber\\
  \equiv &\; \frac{f_1f_6^2}{f_2f_3}\left(
  \sum_{n=0}^\infty
   (-1)^{n(n+1)/2}q^{3n^2+2n}
  \right)^2+\frac{f_1^2f_6^4}{f_2^2f_3^2}
    \sum_{n=0}^\infty(-1)^{n(n+1)/2}q^{3n^2+2n}
  \nonumber\\
  \equiv &\; \frac{f_1f_6^2}{f_2f_3}
   \sum_{n=0}^\infty q^{6n^2+4n}
    +\frac{f_1^2f_6^4}{f_2^2f_3^2}
    \sum_{n=0}^\infty q^{3n^2+2n} \pmod 2.
    \label{2-58}
\end{align}
In view of \eqref{2-12}, \eqref{2-51},
 \eqref{2-57}
  and \eqref{2-58},
\begin{align}
\sum_{n=0}^\infty
 b_3(n)q^n\equiv
 \frac{f_3}{f_1}\cdot \frac{f_4}{f_6f_{24}}
 \sum_{n=0}^\infty q^{6n^2+4n}
  +\frac{f_2}{f_6f_{12}}\sum_{n=0}^\infty
   q^{3n^2+2n} \pmod 2.\label{2-59}
\end{align}
By \eqref{2-12} and \eqref{2-13},
\begin{align}\label{2-60}
\frac{f_3 }{f_1}\equiv \frac{ f_8
   }{f_6}+q\frac{f_{6} f_{24}
  }{f_4 } \pmod 2.
 \end{align}
 Substituting \eqref{a-3}
  and \eqref{2-60}
  into \eqref{2-59} and
  extracting those
terms in which the power of  $q
 $ is congruent to 1 modulo 2, then
 dividing  them by  $q$  and replacing
$q^{2}$ by $q$ yields
\begin{align}\label{2-61}
\sum_{n=0}^\infty
 b_3(2n+1)q^n\equiv
 \sum_{n=0}^\infty q^{3n^2+2n}
  +q^2\cdot\frac{f_3}{f_1}
  \cdot
   \frac{f_2}{f_{12}}\sum_{n=0}^\infty
   q^{6n^2+8n} \pmod 2.
   \qquad ({\rm by}\ \eqref{2-12})
\end{align}
If we substitute \eqref{a-3}
 and  \eqref{2-60}
  into \eqref{2-61} and
  extract  those
terms in which the power of  $q
 $ is congruent to 1 modulo 2, then
 divide  them by  $q$  and replace
$q^{2}$ by $q$, we arrive at
\begin{align}
\sum_{n=0}^\infty
 b_3(4n+3)q^n\equiv
 q^2 \sum_{n=0}^\infty q^{6n^2+8n}
  +q \cdot\frac{f_3^3}{f_1} \cdot
    \sum_{n=0}^\infty
   q^{3n^2+4n} \pmod 2.\qquad ({\rm by}\ \eqref{2-12}) \label{2-62}
\end{align}
Substituting \eqref{2-13} and \eqref{2-40}
  into \eqref{2-62}   yields
\begin{align}\label{2-63}
\sum_{n=0}^\infty
 b_3(4n+3)q^n\equiv &\;
  q  \frac{f_4^3f_6^2}{f_2^2f_{12}}
  \sum_{n=0}^\infty q^{12n^2+8n
   }+q^2 \sum_{n=0}^\infty q^{6n^2+8n}
  +q^2\frac{f_{12}^3
    }{f_4}\sum_{n=0}^\infty q^{12n^2+8n
   }\nonumber\\
   &\;+q^8\frac{f_4^3f_6^2}{f_2^2f_{12}}
  \sum_{n=0}^\infty q^{12n^2+20n
   }+q^9\frac{f_{12}^3
    }{f_4}\sum_{n=0}^\infty q^{12n^2+20n
   }
    \pmod 2.
\end{align}
By \eqref{2-12}
 and \eqref{2-63},
\begin{align}
\sum_{n=0}^\infty
 b_3(8n+3)q^n\equiv
  &\;
  q  \sum_{n=0}^\infty q^{3n^2+4n}
  +q\frac{f_{6}f_{12}
    }{f_2}\sum_{n=0}^\infty q^{6n^2+4n
   } +q^4f_4
  \sum_{n=0}^\infty q^{6n^2+10n
   }
  \pmod 2  \label{2-64}
\end{align}
and
\begin{align}\label{2-65}
\sum_{n=0}^\infty
 b_3(8n+7)q^n\equiv
 f_4 \sum_{n=0}^\infty q^{6n^2+4n}
  +q^4
  \frac{f_6^3}{f_2} \sum_{n=0}^\infty
   q^{6n^2+10n} \pmod 2.
\end{align}
Congruence  \eqref{2-52}
 follows from \eqref{2-65}.

Substituting \eqref{2-40}
 into \eqref{2-64}
  and
   extracting those
terms in which the power of  $q
 $ is congruent to 1 modulo 2, then
 dividing  them by  $q$  and replacing
$q^{2}$ by $q$, we arrive at
\begin{align}\label{2-66}
\sum_{n=0}^\infty
 b_3(16n+11)q^n\equiv
 \sum_{n=0}^\infty
   q^{6n^2+4n}+ \frac{f_3}{f_1} \cdot
    f_6
   \sum_{n=0}^\infty q^{3n^2+2n}  \pmod 2.
\end{align}
Substituting \eqref{a-3}
 and \eqref{2-60} into \eqref{2-66} yields
\begin{align*}
\sum_{n=0}^\infty
 b_3(16n+11)q^n\equiv&\;
 \sum_{n=0}^\infty
   q^{6n^2+4n}+ f_8\sum_{n=0}^\infty
    q^{12n^2+4n}+q \frac{f_6^2f_{24}}{
    f_4}\sum_{n=0}^\infty
    q^{12n^2+4n}
    \nonumber\\
    &\;+q^5f_8
    \sum_{n=0}^\infty
    q^{12n^2+16n} +q^6\frac{f_6^2f_{24}}{
    f_4}\sum_{n=0}^\infty
    q^{12n^2+16n} \pmod 2,
\end{align*}
from which with \eqref{2-12},
 we arrive at
\begin{align}
\sum_{n=0}^\infty
 b_3(32n+11)q^n\equiv
 \sum_{n=0}^\infty
   q^{3n^2+2n} +f_4\sum_{n=0}^\infty
    q^{6n^2+2n}
      +q^3\frac{f_6f_{12}}{
    f_2}\sum_{n=0}^\infty
    q^{6n^2+8n}   \pmod 2 \label{2-69}
\end{align}
and
\begin{align}\label{2-70}
\sum_{n=0}^\infty
 b_3(32n+27)q^n\equiv
         \frac{f_6f_{12}}{
    f_2}\sum_{n=0}^\infty
    q^{6n^2+2n}
    +q^2f_4
    \sum_{n=0}^\infty
    q^{6n^2+8n}   \pmod 2.
\end{align}
Congruence \eqref{2-53}
 follows from \eqref{2-70}.

If we substitute \eqref{a-3}
into \eqref{2-69} and
  extract  those
terms in which the power of  $q
 $ is congruent to 1 modulo 2, then
 divide  them by  $q$  and replace
$q^{2}$ by $q$, we arrive at
\begin{align}
\sum_{n=0}^\infty
 b_3(64n+43)q^n\equiv
  q^2 \sum_{n=0}^\infty
   q^{6n^2+8n} + q  \frac{f_3 }{
    f_1}\cdot f_6\sum_{n=0}^\infty
    q^{3n^2+4n}
     \pmod 2. \label{2-71}
\end{align}
Substituting \eqref{2-40}
 and \eqref{2-60} into \eqref{2-71}
 and
   extracting those
terms in which the power of  $q
 $ is congruent to 1 modulo 2, then
 dividing  them by  $q$  and replacing
$q^{2}$ by $q$, we arrive at
\begin{align}\label{2-72}
\sum_{n=0}^\infty
 b_3(128n+107)q^n\equiv
 f_4\sum_{n=0}^\infty q^{6n^2+4n}
 +q^4\frac{f_3^2f_{12}}{f_2}\sum_{n=0}^\infty q^{6n^2+10n}  \pmod 2.
\end{align}
Congruence \eqref{2-54} follows from
 \eqref{2-12}
  and \eqref{2-72}. This completes
   the proof of this lemma. \qed

\begin{lemma}\label{L-6}
 For $n\geq 0$,
 \begin{align}
b_1(196n+r) &\;\equiv 0 \pmod 4,\label{a-4}
 \end{align}
where $r\in \{110,138,194\}$.

\end{lemma}

\noindent{\it Proof.}
 It follows from \eqref{2-11-0}
 that
 \begin{align}
\sum_{n=0}^\infty
 b_1(2n)q^n =\frac{f_3^3}{f_1}
 \cdot \frac{f_2f_4f_{12}^2}{
 f_6^7}.\label{a-5}
 \end{align}
Substituting \eqref{2-13}
  into \eqref{a-5} and
  extracting those
terms in which the power of  $q
 $ is congruent to 1 modulo 2, then
 dividing  them by  $q$  and replacing
$q^{2}$ by $q$ yields
\begin{align}\label{a-7}
\sum_{n=0}^\infty
 b_1(4n+2)q^n =  \frac{f_1f_{6}^5}{
 f_3^7}.
\end{align}
In view of  \eqref{2-12} and \eqref{a-7},
\begin{align}\label{a-8}
b_1(4n+2)\equiv b(n) \pmod 4,
\end{align}
where $b(n)$ is defined by
\[
\sum_{n=0}^\infty
 b(n)q^n:=f_1f_3f_6.
\]
Utilizing
     the Mathematica package
     RaduRK of   Smoot \cite{Smoot},
 one can prove that for $n\geq 0$,
\begin{align}\label{a-9}
b(49n+r)
 \equiv 0 \pmod 4,
\end{align}
where $r\in\{27,34,48\}$. Replacing
 $n$ by $49n+r $  $(r\in\{27,34,48\})$
  in \eqref{a-8}
   and using \eqref{a-9}, we get \eqref{a-4}.
    This completes the proof
     of Lemma \ref{L-6}. \qed

\begin{lemma}\label{L-4}
 For $n\geq 0$,
\begin{align}\label{2-73}
b_2(196n+s)\equiv
 0 \pmod 2,
\end{align}
where $s\in\{26,54,110,138,166,194\}$.
\end{lemma}

\noindent{\it Proof.} In view of \eqref{2-0},
 \eqref{2-12}, \eqref{2-26}
  and \eqref{2-57},
\begin{align*}
\sum_{n=0}^\infty b_2(n)q^n \equiv &\;
 \frac{ f_2^2f_3^2
 }{f_1^2f_6^4}
  \sum_{n=1}^\infty (-1)^{n(n-1)/2} q^{3n^2-2n}
  \nonumber\\
  \equiv &\; \frac{ f_4f_6
 }{f_2f_{24}}
 \left( \sum_{n=1}^\infty   q^{12n^2-4n}+q
  \sum_{n=0}^\infty   q^{12n^2+8n}\right)
\pmod 2,
\end{align*}
which yields
\begin{align}\label{2-74}
\sum_{n=0}^\infty b_2(2n)q^n
  \equiv \frac{f_3}{f_1}\cdot
   \frac{ f_2 f_6^2
 }{ f_{24}}
  \sum_{n=1}^\infty   q^{6n^2-2n}
\pmod 2. \qquad ({\rm by}\ \eqref{2-12})
\end{align}
Substituting \eqref{2-60} into \eqref{2-74}
 and
   extracting those
terms in which the power of  $q
 $ is congruent to 1 modulo 2, then
 dividing  them by  $q$  and replacing
$q^{2}$ by $q$, we arrive at
\begin{align}
\sum_{n=0}^\infty b_2(4n+2)q^n
  \equiv \frac{f_3^3}{f_1}
  \sum_{n=1}^\infty   q^{3n^2-n}
\pmod 2. \qquad ({\rm by}\ \eqref{2-12}) \label{2-75}
\end{align}
In view of \eqref{2-12},
 \eqref{2-55}
 and \eqref{2-57},
\begin{align}
  \frac{f_3^3}{f_1}
 \equiv \sum_{n=-\infty}^\infty q^{3n^2+2n}
 \pmod 2. \label{2-76}
\end{align}
Combining \eqref{2-75} and \eqref{2-76}
 yields
 \begin{align}
\sum_{n=0}^\infty b_2(4n+2)q^n
  \equiv  \sum_{n=-\infty}^\infty q^{3n^2+2n}
  \sum_{n=1}^\infty   q^{3n^2-n}
\pmod 2. \label{2-77}
\end{align}
It is easy to check that
\begin{align}
\sum_{n=-\infty}^\infty q^{3n^2+2n} =&\;\sum_{n=-\infty}^\infty
 q^{147n^2+14n}
  +\sum_{n=-\infty}^\infty
 q^{147n^2-28n+1} +\sum_{n=-\infty}^\infty
 q^{147n^2+56n+5} +\sum_{n=-\infty}^\infty
 q^{147n^2-70n+8}\nonumber\\
 &\;  +\sum_{n=-\infty}^\infty
 q^{147n^2+98n+16}
 +\sum_{n=-\infty}^\infty
 q^{147n^2-112n+21} +\sum_{n=-\infty}^\infty
 q^{147n^2+140n+33} \label{2-78}
\end{align}
and
\begin{align}
\sum_{n=1}^\infty q^{3n^2-n} =&\; \sum_{n=1}^\infty
 q^{147n^2-7n}
  +\sum_{n=0}^\infty
 q^{147n^2+35n+2} +\sum_{n=0}^\infty
 q^{147n^2+77n+10} +\sum_{n=0}^\infty
 q^{147n^2+119n+24}\nonumber\\
 &\;  +\sum_{n=0}^\infty
 q^{147n^2+161n+44}
 +\sum_{n=0}^\infty
 q^{147n^2+203n+70} +\sum_{n=0}^\infty
 q^{147n^2+245n+102}.\label{2-79}
\end{align}
 Substituting \eqref{2-78} and  \eqref{2-79}
 into \eqref{2-77} and
   extracting those
terms in which the power of  $q
 $ is congruent to 6 modulo 7, then
 dividing  them by  $q^6$  and replacing
$q^{7}$ by $q$, we obtain
\begin{align*}
\sum_{n=0}^\infty b_2(28n+26)q^n \equiv
 \sum_{n=-\infty}^\infty\sum_{m=0}^\infty
 q^{21n^2+14n+21m^2+35m+16},
\end{align*}
which yields \eqref{2-73}. This completes
 the proof of this lemma. \qed

Now, we are ready to prove Theorem
 \ref{Th-1}.

\noindent{\it Proof
 of Theorem \ref{Th-1}.}
  Note that
\begin{align}
\frac{1}{A-B}&\;=\frac{A-B}{(A+B)^2-4AB}
 \nonumber\\
&\;=\frac{(A-B)((A+B)^2+4AB)}{(A+B)^4-16A^2B^2}
\nonumber\\
&\;\equiv \frac{(A+B-2B)((A+B)^2+4AB)}{(A+B)^4}
 \nonumber
\\
&\;\equiv \frac{1}{A+B}
 -\frac{2B}{(A+B)^2}
 +4\frac{A^2B+AB^2}{(A+B)^4}
\pmod 8, \label{3-1}
\end{align}
where $A$ and $B$ are defined by \eqref{2-0}.

In view of \eqref{1-6},
 \eqref{1-7}, \eqref{2-0} and \eqref{3-1},
 \begin{align*}
\sum_{n=0}^\infty c_5(n)q^n\equiv
 \frac{1}{A+B}
 -\frac{2B}{(A+B)^2}
 +4\frac{A^2B+AB^2}{(A+B)^4}
\pmod 8,
 \end{align*}
from which with \eqref{2-1}, \eqref{2-26}
 and \eqref{2-51},
 we deduce that for $n\geq 0$,
 \begin{align}\label{3-2}
c_5(n)\equiv b_1(n)-2b_2(n)+4b_3(n) \pmod 8.
 \end{align}
Congruences \eqref{1-8}--\eqref{1-10} follow
 from \eqref{3-2}
  and Lemmas \ref{L-1}--\ref{L-3}.

 In addition, it follows from \eqref{3-2} that
 \begin{align}\label{3-3}
c_5(n)\equiv b_1(n)-2b_2(n) \pmod 4.
 \end{align}
 Congruences \eqref{1-11} and \eqref{1-12} follow
 from \eqref{3-3},
  and Lemmas \ref{L-1} and \ref{L-2}.

  Furthermore, Congruence
   \eqref{1-13} follows
 from  \eqref{a-4}, \eqref{2-73}
  and \eqref{3-3}. The proof of Theorem
   \ref{Th-1}
    is complete. \qed

\section{Proof of Theorem \ref{Th-2}}

To prove Theorem \ref{Th-2},
 we first prove two lemmas.

\begin{lemma}\label{L-7}
If $p\equiv 7\pmod 8 $ is a prime, then
 for all   integers
  $n$ and $\alpha\geq 0$ with $p\nmid n$,
\begin{align}\label{4-1}
b_1\left(4
 p^{2\alpha+1}n+\frac{8 p^{2\alpha+2}+1
  }{3}
 \right)\equiv 0 \pmod
4.
\end{align}

\end{lemma}

\noindent{\it Proof.} By \eqref{2-12} and \eqref{2-15},
 \begin{align*}
\sum_{n=0}^\infty b_1(4n+3)q^n\equiv&\;
 -  \frac{f_3^6}{f_1^2}\cdot\frac{f_3^2}{f_6
  } \nonumber\\
  \equiv &\;
 -\sum_{n=0}^\infty
 a_1(n)q^n +2\sum_{n=0}^\infty
  a_2(n)q^n \pmod 4, \qquad
  ({\rm by}\ \eqref{1-4})
  \end{align*}
  where
  \begin{align}
\sum_{n=0}^\infty
 a_1(n)q^n:=&\;
  \frac{f_3^6}{f_1^2},\nonumber\\
  \sum_{n=0}^\infty a_2(n)q^n:=&\;\frac{f_6^3}{f_2}
  \sum_{n=1}^\infty
   q^{3n^2}. \label{v-2}
  \end{align}
Thus,
 for $n\geq 0$,
 \begin{align}\label{4-2}
b_1(4n+3)\equiv -a_1(n)+2a_2(n) \pmod 4.
 \end{align}
In \cite{Wang-0}, Wang  proved that if   we write
 $3n+2=\prod_{i=1}^sp_i^{\alpha_i}
 $ as the
unique prime factorization, then
 \begin{align}\label{4-3}
a_1(n)=\frac{1}{3}
 \prod_{i=1}^s \frac{p_i^{\alpha_i+1}-1}{ p_i-1}.
\end{align}
In view of \eqref{4-3}, we see that if $
 p\equiv 7\pmod 8 $ is a prime, then
 for all   integers
  $n$ and $\alpha\geq 0$ with $p\nmid n$,
\begin{align}\label{4-4}
a_1\left(
 p^{2\alpha+1}n+\frac{2(p^{2\alpha+2}-1)}{3}\right)
  \equiv 0 \pmod 4.
\end{align}
Thanks to \eqref{2-76} and \eqref{v-2},
\begin{align*}
  \sum_{n=0}^\infty a_2(n)q^n\equiv
   \sum_{n=-\infty}^\infty q^{6n^2+4n}
  \sum_{n=1}^\infty
   q^{3n^2} \pmod 2,
\end{align*}
which yields
\[
  \sum_{n=0}^\infty a_2(n)q^{3n+2}\equiv
   \sum_{n=-\infty}^\infty \sum_{m=1}^\infty
     q^{2(3n+1)^2+(3m)^2} \pmod 2.
\]
Therefore, if $3n+2$ is not of
 the form $2x^2+y^2$,
  then
\[
 a_2(n)\equiv 0 \pmod 2.
\]
Note that if $N$ is   of the form $2x^2+ y^2$,
 then $\nu_p(N)$ is even since
    $p$ is a prime
 with $p\equiv 7 \pmod 8 $ and  $\left(
 \frac{-2 }{p}\right)=-1$.  Here and throughout
  this paper, $\nu_p(N)$
  denotes the highest power of $p$
   dividing
$N$ and $\left(
 \frac{ \cdot }{p}\right)$ denotes the Legendre symbol.
  It is easy to check that
   if $p\nmid n$, then
\[
 \nu_p\left( 3\left(
 p^{2\alpha+1}n+\frac{2(p^{2\alpha+2}-1
  )}{3}\right)+2\right)
  =\nu_p\left(    p^{2\alpha+1}(3n+2p)
  \right)
  =2\alpha+1
\]
is odd. Therefore,
 $3\left(p^{2\alpha+1}n+\frac{2(p^{2\alpha+2}-1
  )}{3}\right)+2\ (p\nmid n)$ is not of the form
  $2x^2+ y^2$ and
\begin{align}\label{4-7}
a_2\left( p^{2\alpha+1}n+\frac{2(p^{2\alpha+2}-1
  )}{3}\right)\equiv 0 \pmod 2.
\end{align}
Replacing $n$ by $ p^{2\alpha+1}n+\frac{2(p^{2\alpha+2}-1
  )}{3}\ (p\nmid n)$
   in \eqref{4-2} and utilizing \eqref{4-4}
     and \eqref{4-7}, we arrive
   at \eqref{4-1}. This completes
    the proof.
 \qed

\begin{lemma} \label{L-5}
If $p\equiv 7\pmod 8 $ is a prime, then
 for all   integers
  $n$ and $\alpha\geq 0$ with $p\nmid n$,
\begin{align}\label{2-80}
b_2\left(4
 p^{2\alpha+1}n+\frac{8p^{2\alpha+2}+1
  }{3}\right) \equiv
0 \pmod 2.
\end{align}
\end{lemma}

\noindent{\it Proof.} By \eqref{2-35} and \eqref{a-3},
\begin{align*}
\sum_{n=0}^\infty
 b_2(4n+3)q^n
&\;\equiv  \frac{f_3^3}{f_1}
 \cdot   \sum_{n=0}^\infty
  q^{3n^2+2n}
   \pmod 2 , \qquad ({\rm by}\ \eqref{2-12})
\end{align*}
from which with \eqref{2-76},
 we get
\begin{align*}
\sum_{n=0}^\infty
 b_2(4n+3)q^n
&\;\equiv \sum_{n=-\infty}^\infty q^{3n^2+2n}  \sum_{n=0}^\infty
  q^{3n^2+2n}
   \pmod 2.
\end{align*}
From the above congruence,
 we arrive at
 \begin{align*}
\sum_{n=0}^\infty
 b_2(4n+3)q^{3n+2}
&\;\equiv \sum_{n=-\infty}^\infty \sum_{m=0}^\infty
  q^{(3n+1)^2+(3m+1)^2}
 \pmod 2.
\end{align*}
Therefore, if $3n+2$ is not of
 the form $x^2+y^2$,
  then
\[
 b_2(4n+3)\equiv 0 \pmod 2.
\]
Note that if $N$ is   of the form $x^2+ y^2$,
 then $\nu_p(N)$ is even since
    $p$ is a prime
 with $p\equiv 7 \pmod 8 $ and  $\left(
 \frac{-1 }{p}\right)=-1$.   It is easy to check that
   if $p\nmid n$, then
\[
 \nu_p\left( 3\left(
 p^{2\alpha+1}n+\frac{2(p^{2\alpha+2}-1
  )}{3}\right)+2\right)
  =\nu_p\left(    p^{2\alpha+1}(3n+2p)
  \right)
  =2\alpha+1
\]
is odd. Therefore,
 $3\left(p^{2\alpha+1}n+\frac{2(p^{2\alpha+2}-1
  )}{3}\right)+2\ (p\nmid n)$ is not of the form
  $x^2+ y^2$ and
  \eqref{2-80} holds. The proof is complete.
   \qed

To conclude this section, we pose
 a proof of Theorem \ref{Th-2}.

 \noindent{\it Proof of Theorem
 \ref{Th-2}.} Replacing $n$
  by $4
 p^{2\alpha+1}n+\frac{8 p^{2\alpha+2}+1
  }{3}$ in \eqref{3-3}
   and using \eqref{4-1}
    and \eqref{2-80}, we arrive at \eqref{1-15}.
     The proof of Theorem \ref{Th-2}
      is complete. \qed

\section{Concluding remarks}

 As seen in Introduction, false theta functions
   have received
 a lot of attention
 in recent years.
 In this study, we confirm a conjecture
 of Keith on congruences
  of modulo 4 and 8   for the coefficients
        of the  reciprocal  of a
         false theta function. A natural
   question
 is to extend
    the congruences in
     this paper  to modulo $16$, $32$,
      $64$, etc.
 However,  it will likely
  require a different approach
  since the methods used
in this paper run into serious
 limitations beyond the modulus
   of $16$. Moreover,  it would be
    interesting
     to determine the
      arithmetic density of   the
        set of
       integers such that
        $ c_t(n)
        \equiv 0 \pmod
         {m}$ for   fixed positive
          integers $t$  and   $m\geq 2$.

\section*{Statements and Declarations}

 \noindent{\bf Funding.}
   This work was supported by
 the National Natural Science Foundation of
  China  (grant no.
    12371334) and the Natural Science Foundation of
   Jiangsu Province of China (grant no.
    BK20221383).

\noindent{\bf Author Contributions.}
      The authors contributed equally to the preparation of this article. All authors read and
 approved the final manuscript.

 \noindent{\bf Competing Interests.}
 The authors declare that they have
no conflict of interest.

  \noindent{\bf Data Availability Statements.} Data sharing not applicable to this
article as no datasets were generated or analysed during the current
study.

\end{document}